\input amstex
\documentstyle{amsppt}
\magnification \magstep 1 
\def\Zt{\Bbb Z}
\def\Pt{\Bbb P}
\def\Rt{\Bbb R}
\def\Ct{\Bbb C}

\def\Hc{\Cal H}
\def\l{\ell}
\def\iii{}
\def\D#1{\,{\hbox{\it D\kern-9pt I\kern5.5pt}}^{_{^\bullet}}_{#1}}
\def\m{\setminus}
\def\Fc{\Cal H^{_{^\bullet}}}
\def\cF{\Cal F}
\def\cP{\Cal P_{\Cal V}}
\def\codim#1#2{\mathord{\mathop{\roman{codim(}}\nolimits}#1\smash\subset\smash\, #2)}
\def\Hom{\mathord{\mathop{\roman{Hom}}\nolimits}}
\def\HH{{\hbox{\it H\kern-9pt I\kern6.5pt}}}
\def\map{\smash{\mathop{\hbox to 25pt{\rightarrowfill}}}} 
\def\mono#1{\mathop{\hookrightarrow}\limits^#1}
\def\epi#1{{\mathop{\hbox to
13pt{\rightarrowfill}}\limits^{#1}}\hskip -11pt\rightarrow} 
\def\i#1{\Cal{IC}^{_{^\bullet}}_{#1}} 
\def\p#1{\smash{\Cal I\hskip-1pt_{\goth p}  \Cal C^{_{^\bullet}}_{#1}} }
\def\q#1{\smash{\Cal I_{\underline  k} \Cal C^{_{^\bullet}}_{#1}}}
\def\h#1#2{IH_{#1}({#2})}
\def\hc#1{IH^\roman{cld}_{#1}}

\topmatter
\title A Morphism of Intersection Homology \\ and Hard Lefschetz 
\endtitle 
\author Andrzej Weber \endauthor
\thanks Partially supported by KBN 2 P03A 01113 grant\endthanks

\affil Institute of Mathematics, University of Warsaw \endaffil
\address ul. Banacha 2, 02--097 Warszawa, Poland \endaddress
\email aweber\@mimuw.edu.pl \endemail

\abstract We consider a possibility of the existence of intersection
homology morphism, which would be associated to a map of analytic
varieties. We assume that the map is an inclusion of codimension one.
Then the existence of a morphism follows from Saito's decomposition
theorem.  For varieties with conical singularities we show, that the
existence of intersection homology morphism is exactly equivalent to
the validity of Hard Lefschetz Theorem for links. For varieties with
arbitrary analytic singularities we extract a remarkable property,
which we call Local Hard Lefschetz. 
\endabstract

\endtopmatter

\document 
\head 0. Introduction\endhead
Let $Y$ be a complex algebraic variety of pure dimension.  Any
algebraic subvariety $X$ 
of dimension $i$ defines a class in the homology group with closed
supports: $$[X]\in H^\roman{cld}_{2i}(Y;\Zt)\,.$$
It was conjectured in \cite{Bry}, \S5 Probl\`emes ouverts, that
any such class lifts to
intersection homology of $Y$. There are examples that it is not the
case for integer coefficients. For rational coefficients the answer
is positive. Barthel et al.~(\cite{BBFGK}) proved that for any morphism
of algebraic varieties (not only an inclusion) there exists an
associated homomorphism of rational intersection homology, however this morphism 
is not unique.
Unfortunately the methods of the existing proofs are extremely
advanced. The key argument is to reduce to a finite characteristic and
then to apply the results of \cite{BBD}.  
Another proof is made by means of resolution of singularities, see
\cite{We}. It relies on the decomposition theorem, which again
follows from \cite{BBD}.

Any inclusion of algebraic varieties may be factorized by 
inclusions of codimension one. For general {\it analytic} varieties 
we cannot make this
reduction, but we will assume that we are
given a pair of analytic varieties with $\codim X Y=1$.
We will only consider intersection homology with rational coefficients.
This paper has three purposes. First we want to show how to
deduce the existence of intersection homology morphisms from the Hard
Lefschetz Theorem provided that singularities are conical, see
Definition 1.1 and Theorem 1.2. Our method is related to the
earlier paper of Brasselet, Fieseler and Kaup \cite{BFK}. 
Later, we give a proof for general analytic singularities
($\codim X Y=1$), which is based on the decomposition theorem proved
by M.~Saito for projective morphisms. 

The second purpose is the following~: Proposition 2.2 asserts that an
intersection homology morphism extends to a stratum if and only if the map of
intersection homology of its links vanishes in the middle dimension.
This is always the case for analytic varieties.  \iii
The vanishing property are stated in Corollaries
5.1--2. For isolated singularities we have:

\proclaim{Proposition 0.1} Let $X\subset Y$ be a pair of analytic
varieties with $\dim X=n$ and $\dim Y=n+1$. Suppose that $x$ is an
isolated singular point of $X$ and $Y$. Let $\l_X\subset\l_Y$ be its
links in $X$ and $Y$. Then the maps induced by the inclusion of links
$\iota_*:H_n(\l_X)@>>>H_n (\l_Y)$ vanishes.
\endproclaim

This assertion can also be derived from the mixed Hodge structure on the 
links.
For nonisolated singularities, it holds when one replaces homology by
intersection homology. Moreover, this property of links implies the
Hard Lefschetz Theorem as shown in Proposition 5.4.  We are tempted to
call this property ''local Hard Lefschetz''. In the algebraic case it
follows directly from Lemme clef of \cite{BBFGK}, 3.3.

Finally we want to show a new method of constructing elements in
intersection homology. Any sequence of subvarieties
$$X=X^0\supset X^1\supset \dots \supset X^k$$ 
with $\codim {X^i} X=i$ such that no component
of $X^i$ is 
contained in the singularities of $X^{i-1}$ determines a sequence of
classes in intersection homology. This method may be applied to lift
Chern classes to intersection homology as it is done in \cite{BW}.

Let us comment the argument for conical singularities.
Our proof of the existence of intersection homology morphisms is dual to
the one in \cite{BBFGK}. We expose the geometry which is hidden behind
the operation on sheaves. We also state a result about uniqueness.
The proof is by induction with respect to strata. The inductive step has
two parts. 
The first: Propositions 2.1 and 2.2 follows from standard operation
on sheaves.
It remains true for any real stratified space with even dimensional strata. 
The geometry is explained in Remark 2.5. The mysterious part is Proposition 
3.1 which makes sense only for conical singularities. 
Our approach is a step towards a full geometric proof, but such a proof
cannot exists unless there is a geometric proof of the Hard Lefschetz.

I would like to thank Professor Gottfried Barthel for his patience, valuable
comments and various improvements of this paper.

\head 1. Morphism of Intersection Homology\endhead

Let $(Y,X)$ be a pair of complex analytic varieties with
$\codim X Y=1$. Let $\i Y$ be the intersection
homology sheaf with rational coefficients (which is an object 
of the derived category of sheaves on $Y$). We follow the 
convention according to which $\i Y$ restricted to the nonsingular part 
of $Y$ 
is concentrated in the dimension zero. The coefficients are rationals.
Then intersection homology groups with closed supports (for the middle
perversity $\goth m$, which is always tacitly understood unless
otherwise stated) are hypercohomology groups \cite{GM1}:
$$\hc i(Y) = \HH^{2\dim Y-i}(Y;\i Y)\,.$$
Let $\i X$ be the intersection homology  sheaf of $X$ considered as a 
sheaf on $Y$ supported by $X$. Then
$$\hc i(X) = \HH^{2\dim X-i}(Y;\i X)\,.$$
Analogously let $\D Y=C^\roman{cld}_{2\dim Y - {_{^\bullet}}}(Y)$ and $\D
X=C^\roman{cld}_{2\dim X - {_{^\bullet}}}(X)$ be the sheaves of Borel-Moore
chains. The homology map is induced by a sheaf morphism:
$$\roman{incl}_\#:\D X[-2] @>>> \D Y\,,$$
where $$\left(\D X[-2]\right)_i=\left(\D X\right)_{i-2}=C^\roman{cld}_{2\dim
X - i+2}(X)=C^\roman{cld}_{2\dim Y - i}(X)\,.$$
We are looking for a morphism 
$$\alpha_{X,Y}: \i X[-2]@>>> \i Y\,,$$
which is {\it compatible with homology}, i.e.~makes the
following diagram commutative:
$$\CD \i X[-2] @>\alpha_{X,Y}>> \i Y\\
@VVV @VVV\\
\D X[-2] @>\roman{incl}_\#>> \D Y\,.
\endCD$$
The vertical arrows in the diagram are the canonical morphisms from 
intersection homology  to homology.
Assume that no component of $X$ is contained in the singularities of $Y$. 
Let $\p Y$ be the intersection homology sheaf for the logarithmic
perversity $\goth p$:
$${\goth p}(2c)={\goth m}(2c)+1=c\qquad\roman{for}\quad c>0\,.$$
Suppose that $\Cal S$ is a Whitney stratification of the pair
$(Y,X)$. Let $Y_\roman{reg}\,$ be the nonsingular part of $Y$. 
Construct another stratification $\Cal S'$ of $Y$ for which
$Y_\roman{reg}\,$ is the biggest stratum and $Y\m Y_\roman{reg}\,$
is stratified by $\Cal S$.  If a geometric chain in $X$ is allowable with
respect to $\goth m$ and to the stratification $\Cal S\cap X$, then
it is allowable in $Y$  with respect to $\goth p$ and $\Cal S'$.
Thus the inclusion of geometric chains $C_{_{^\bullet}}^\roman{cld}(X)\subset
C_{_{^\bullet}}^\roman{cld}(Y)$ gives us a morphism of sheaves 
$$\roman{Incl}_\#:\i X[-2] @>>>\p Y\,.$$ 
We will say that a morphism $\alpha_{X,Y}$ is {\it geometric} if
the following triangle commutes:
$$\matrix  \i X[-2] & @>\alpha_{X,Y}>> & \i Y\\
\hfill_{\roman{Incl}_\#}& \searrow\hfill& @VVV\\
 & & {\p Y}^{\phantom\int}.&
\endmatrix$$
The vertical arrow is the one induced by the inequality $\goth m < \goth p$.
If $\alpha_{X,Y}$ is geometric, then it is compatible with homology.
On $Y_\roman{reg}\,$ we have a
geometric morphism $\alpha_\roman{reg}:\i{X\cap
Y_\roman{reg}}[-2]@>>> \i{ Y_\roman{reg}}$ 
which is the composition of the natural morphisms:
$$\i{X\cap Y_\roman{reg}}[-2]@>>>\D{X\cap Y_\roman{reg}}[-2]
@>\roman{incl}_\#>>\D{Y_\roman{reg}}\cong\i{ Y_\roman{reg}}\,.$$ 
This is the unique morphism which is compatible with homology by construction.
We are looking for an extension of $\alpha_\roman{reg}\,$ to all $Y$. 
We will construct it for a case where the singularities have a special form.

\remark{Definition 1.1} Let  $X\subset  Y\subset  M$  be  a  pair  of 
analytic
subvarieties of a smooth analytic manifold. We say that it has {\it 
conical singularities} if there exists a Whitney stratification
of $(Y,X)$ and for (each connected component of) each stratum $S$, there 
exists
a normal slice $N_S$ in $M$ with coordinates such that $N_S\cap X$
and $N_S\cap Y$ are both given by homogeneous equations. \endremark

\iii
Examples of such spaces are Schubert varieties in Grassmannians.
From the Hard Lefschetz of \cite{BBD} we will deduce the following:

\proclaim{Theorem 1.2} Let $X\subset Y$ be a pair of analytic varieties
with conical singularities and $\codim X Y=1$.
Suppose that no 
component of $X$ is contained in the singularities of $\;Y$. Then there
exists an unique 
extension $\alpha_{X,Y}$ of $\alpha_\roman{reg}\,$ and this morphism is
geometric.\endproclaim 

The case of general analytic singularities follows from M.~Saito's
theory of mixed Hodge modules. We give both proofs in \S3.
The unicity is in fact proved between the lines of
\cite{BBFGK}, although it was not explicitly stated.

\remark{Remark 1.3} If a component of $X$ is contained in the
singularities of $Y$, then to construct a morphism of intersection
homology, we would have to go through the
normalization process of \cite{BBFGK}, p.166. At this point we
would loose uniqueness.
\endremark\vskip 5pt

\head 2. Extending an Intersection Homology Morphism\endhead

The proof will be inductive.
For  the inductive step we fix some notation. We assume that the
pair $(Y,X)$ is stratified by a Whitney stratification. Let $\dim X= n$,
$\dim Y=n+1$. 
Let $S$ be a minimal (connected) stratum of $X\m
Y_\roman{reg}\,$. Put $k=\codim S X$. 
Let $X'=X\m S$ and $Y'=Y\m S$  and let $i:
Y'\hookrightarrow Y$ be the inclusion. Then 
$$\i Y=\tau_{\leq k}Ri_*\i {Y'}\,,$$
$$\p Y=\tau_{\leq k+1}Ri_*\p {Y'}\,,$$
and $$\i X=\tau_{\leq k-1}Ri_*\i {X'}\,.$$
Denote by $\q Y$ another intersection homology sheaf on $Y$ which is
associated to the perversity
$${\underline k}(2i)=
\left\{\matrix \goth m(2i)=i-1&\;\roman{for}\;i\leq k\hfill\\
\goth m(2i)+1=i&\;\roman{for}\;i>k.\endmatrix\right.$$
Through the sheaf
$$\q Y=\tau_{\leq k+1}Ri_*\i {Y'}\,,$$
factors the canonical morphism $\i Y@>>>\p Y$. 
Suppose we are given a morphism 
$$\alpha_{X',Y'}:\i {X'}[-2]@>>>\i {Y'}\,.$$ 
Applying the functor $\tau_{\leq k+1}Ri_*$ to $\alpha_{X',Y'}$ we
obtain a morphism $\phi:\i X[-2]@>>>\q Y$ as follows:
$$\matrix \tau_{\leq k+1}Ri_*\alpha_{X',Y'}:& 
\tau_{\leq k+1}Ri_*(\i {X'}[-2]) &@>>>& \tau_{\leq k+1}Ri_*\i {Y'}\\ 
&\parallel& &\parallel\\ 
&(\tau_{\leq k-1}Ri_*\i {X'})[-2]&&\q Y\\
&\parallel&\;\nearrow_\phi\\ 
&\i X[-2]\endmatrix$$

\proclaim{Proposition 2.1} The restriction of sheaves to $Y'$
induces an isomorphism
$$\Hom(\i X[-2],\q Y)\cong \Hom(\i{X'}[-2],\i{Y'})\,.$$
In particular, every extension of $\alpha_{X',Y'}$ equals $\phi$.
\endproclaim 

This proposition follows from \cite{BBFGK}, Lemme 3.1. For the sake
of completeness we will give a proof in our case.

\demo{Proof} We have 
$$\aligned \Hom(\i {X'}[-2],\i {Y'})\cong&
\;\Hom(Ri^*\i X[-2],\i {Y'})\cong\\
\cong&\;\Hom(\i X[-2],Ri_*\i {Y'})\,.\endaligned$$
Consider the distinguished triangle:
$$\def\map{\smash{\mathop{\hbox to
60pt{\rightarrowfill}}}} 
\matrix \q Y &  \map &Ri_*\i {Y'}\,.\\
\hfill_{[+1]}\nwarrow & &\swarrow\hfill\\
 &\tau_{\geq k+2}Ri_*\i {Y'}\endmatrix$$
Applying the functor $\Hom(\i X[-2],\,\cdot\,)$, we obtain a long
exact sequence: 
\vskip 5pt
\noindent $@>>>\Hom(\i X[-2],\tau_{\geq k+2}Ri_*\i {Y'}[-1])@>>>
\Hom(\i X[-2],\q Y)@>>>$
\vskip 5pt
\hfill $ @>>>\Hom(\i X[-2],Ri_*\i {Y'})
@>>>\Hom(\i X[-2],\tau_{\geq k+2}Ri_*\i {Y'})@>>>$
\vskip 5pt \noindent
The cohomology of the sheaf $\i X[-2]$ is concentrated in dimensions
$2, 3,\dots ,k+1$. Thus the first and the fourth $\Hom$--groups
vanish, so in the middle we have an isomorphism.\qed\enddemo

Now we want to lift the morphism $\phi$ to $\i Y$. To this end, we
introduce more notation. Let $N_{S}$ be a normal
slice of $S$ in $M$ at some point $x\in S$. 
Then there exists a neighbourhood $U$ of $x$ such that $U\cap X\simeq
B\times c\l_X $ and $U\cap Y\simeq B\times c\l_Y$, where $B$ is a ball
of the dimension $2(n-k)$ and $c\l_X $ (resp. $c\l_Y$) is the cone over
the link $\l_X $ (resp. $\l_Y$) of $S$ in $X$ (resp. $Y$). The morphism
$\alpha_{X',Y'}$ induces a map for the $k$-th intersection homology
of the links:
$$(\alpha_{\l_X ,\l_Y})_k: \h k {\l_X }@>>>\h k {\l_Y}\,.$$ 

\proclaim{Proposition 2.2} The morphism $\alpha_{X',Y'}$
extends to a morphism $\alpha_{X,Y}$
if and only if the map $(\alpha_{\l_X ,\l_Y})_k$ vanishes. If
$\alpha_{X,Y}$ exists, then 
it is unique. If $\alpha_{X',Y'}$ is geometric, then
$\alpha_{X,Y}$ is geometric. \endproclaim

\demo{Proof} 
We recall from \cite{Bo}, \S1 that for a compact pseudomanifold $L$,
we have
$$ 
\hc i(cL) \cong \left\{\matrix \hc i(cL\m\{*\})\cong\h{i-1} L
&\;\roman{for}\quad i>\frac12\dim_\Rt cL\\ 
\quad 0 \hfill &\;\roman{for}\quad i \leq \frac12\dim_\Rt cL\endmatrix \right.
\tag 2.3$$
Assume that the extension $\alpha_{X,Y}$ exists. Then the map
$$(\alpha_{\l_X ,\l_Y})_k:\h k {\l X}@>>>\h k {\l_Y }$$ 
may be completed to the commutative diagram which is induced by
$\alpha_{X,Y}$ and the restriction from $c\l_Y$ to $c\l_Y\m\{x\}$:
$$\matrix & & \hc {k+1}(c\l_X )&@>\cong>>&\hc {k+1}(c\l_X \m\{x\})&\cong&
\h k {\l_X }\;\\
& &  @VVV   @VVV   @VV(\alpha_{\l_X ,\l_Y})_k V\\
0 &=&\hc {k+1}(c\l_Y)&@>>>&
\hc {k+1}(c\l_Y\m\{x\})&\cong&\h k {\l_Y }\,.\endmatrix\tag 2.4$$
The upper horizontal arrow is isomorphism and $\hc {k+1}(c\l_Y)=0$ by 2.3. 
Thus the map $(\alpha_{\l_X ,\l_Y})_k$ vanishes.

To prove the converse,
consider another distinguished triangle: 
$$\matrix \i Y  & @>>> & \q Y\,, \\
 \hfill_{[+1]}\nwarrow & & \swarrow \hfill\\
 & \Fc
\endmatrix$$
where $\Fc = \Hc^{k+1}Ri_*\i {Y'}[-k-1]$ is the obstruction
sheaf \cite{GM1}, \S5.5. As in the proof of 2.1 above, the triangle
induces a long exact sequence: 
\vskip 5pt\noindent $ @>>>\Hom(\i X[-2],\Fc [-1])@>>>
\Hom(\i X[-2],\i Y)@>>>$
\vskip5pt\hfill
$ @>>> \Hom(\i X[-2],\q Y)@>>>\Hom(\i X[-2],\Fc )@>>>$
\vskip5pt\hskip190pt
$\phi\hskip31pt \mapsto \hskip31pt \roman{res}$
\vskip 5pt\noindent 
The morphism $\phi$ factors through $\i Y$ if and only if the resulting
residue morphism $$\roman{res}:\i X[-2]@>>>\Fc$$ vanishes. 
Note that in this case the lift $\i X[-2] @>>>\i Y$ of $\phi$ is unique
(compare \cite{BBFGK}, p.167)~: This is because the first term in the
sequence vanishes as $\i X[-2]$ is concentrated
in dimensions $2, 3,\dots , k+1$ and the target sheaf only lives in 
dimension $k+2$. Moreover, by Proposition 2.1, any extension of
$\alpha_{X',Y'}$ comes from $\phi$.

Now we apply the functor $\tau_{\geq k+1}$ to the residue morphism and
obtain the commutative square
$$\def\map#1{\;\smash{\mathop{\hbox to 40pt{\rightarrowfill}}\limits^{#1}}\;}
\matrix & & \i X[-2]&\map{\roman{res}}&\Fc\\
& & @VVV @VV\cong V\\
\Hc^{k+1}(\i X[-2])[-k-1]&\cong&\tau_{\geq k+1}(\i
X[-2])&\map{\tau_{\geq k+1}\roman{res}}&\tau_{\geq
k+1}\Fc\,.\endmatrix $$ 
We see that the residue morphism factors through $\Hc^{k+1}(\i
X[-2])[-k-1]$.
The source and the target of $\;\tau_{\geq k+1}\roman{res}\;$ are
sheaves which are concentrated in the dimension $k+1$. They are
supported on $S$ and are locally constant on $S$. It suffices to
examine $\;\tau_{\geq k+1}\roman{res}\;$ at one points of $S$.   
The stalks on $S$ of the sheaves considered here are:
$$\Hc^{k-1}(\i X)_x=\hc {2\dim X-k+1}(U)\cong\hc {k+1}(c\l_X 
)\cong\h k {\l_X }$$ 
and
$$\Hc^{k+1}(Ri_*\i {Y'})_x=\hc {2\dim Y-k-1}(U\m S)\cong
\hc {k+1}(c\l_Y\m \{x\})\cong\h k {\l_Y}\,.$$
The morphism of the stalks is
$$(\alpha_{\l_X ,\l_Y})_k:\h k {\l_X }@>>>\h k {\l_Y}\,.$$ 
When it vanishes, then the residue morphism vanishes as well.

If $\alpha_{X',Y'}$ is geometric then
$\alpha_{X,Y}$ is geometric. This 
is because the sheaf $\q Y$ maps to $\p Y$ and the composition with $\phi$
is the morphism $\roman{Incl}_\#$ as desired.\qed\enddemo

\remark{Remark 2.5} To see some geometry hidden  behind  the  sheaves 
consider
the following situation: $n=k$, $S=\{x\}$, $Y\m\{x\}$ is smooth. We
are given a cycle $\xi$ which is allowable in $X$. Our goal is to
make it allowable in $Y$. We have three cases:
\item {1)} If $\dim\xi\leq n$, then $\xi$ is not allowed to intersect
$\{x\}$ in $X$, so it is allowable in $Y$.
\item {2)} If $\dim\xi\geq n+2$, then $\xi$ is allowed to intersect
$\{x\}$ both in $X$ and $Y$.
\item {3)} A problem occurs if $\dim\xi=n+1$, since then $\xi$ is allowed
to intersect $\{x\}$ in $X$ but not in $Y$. In a neighbourhood of $x$
the cycle $\xi$ is the cone over a cycle $\eta\subset\l_X $.
Suppose that $(\alpha_{\l_X ,\l_Y})_k:\h k {\l_X }@>>>\h k {\l_Y}$
vanishes, i.e.~$\eta$ is a boundary in $\l_Y$, say $\eta=\partial\zeta$. Then 
the
cycle $(\xi\m c\eta)\cup_\eta\zeta$ does not intersect $\{x\}$ and it
is homologous with $\xi$ in $Y$. This is the allowable cycle that we
are looking for.
\endremark\vskip 5pt

A similar consideration for isolated singularities one can find
in \cite{Y}.

\head 3. The Main Theorem\endhead

Let ${\goth X}\subset{\goth Y}\subset\Pt^N$ be a pair of projective
varieties. Denote by $X\subset\Ct^{N+1}$ (resp.~
$Y\subset\Ct^{N+1}$) the affine cones over the ${\goth X}$ (resp.~ $\goth
Y$) and by $\l_X$ (resp.~ $\l_Y$) the links at $0\in\Ct^{N+1}$.
We have a
circle action on $\l_X  $ and $\l_Y $. Let $p:\l_X  @>>>\l_X  /S^1={\goth X}$ 
and
$q:\l_Y @>>>\l_Y /S^1={\goth Y}$ be the quotients.
We have isomorphisms:
$$Rq^*\i {\goth X} \cong \i {\l_X }\quad \roman{and}\quad
Rq^*\i {\goth Y}\cong\i {\l_Y }\,.$$ 

\proclaim{Proposition 3.1} Let ${\goth X}\subset{\goth Y}\subset\Pt^N$ be
a pair of projective varieties with $\dim {\goth X}=k-1$ and $\dim {\goth 
Y}=k$.
Suppose that there is given a morphism of intersection
homology $\alpha_{{\goth X}, {\goth Y}}:\i {{\goth X}}[-2]@>>> \i {{\goth 
Y}}$.
Then the induced map of intersection homology 
$(\alpha_{\l_X,\l_Y})_k: \h k {\l_X }@>>>\h k
{\l_Y}$ vanishes.\endproclaim

\demo{Proof}  
We have a morphism of Gysin sequences of the fibrations $p$ and $q$ induced
by $\alpha_{\goth X, \goth Y}$: 
$$\matrix   \h {k+1} {\goth X} &\mono \Lambda& \h {k-1} {\goth X} &\epi{p^*}& 
\h k
{\l_X  } & @>{p_*}>0> & \h k {\goth X} & @>\Lambda>\cong> & \h {k-2} {\goth 
X} \\
 @VVV  @VVV  @VV(\alpha_{\l_X,\l_Y})_kV  @VVV  \Big\downarrow\\
\h {k+1} {\goth Y} &@>\Lambda>\cong>&\h {k-1} {\goth Y} &@>{q^*}>0>& \h k
{\l_Y } & \mono{{q_*}} & \h k {\goth Y} & \epi \Lambda & \h {k-2} {\goth Y}\,.
\endmatrix$$ 
The map $\Lambda$ is given by intersecting with the Chern class
of the tautological 
bundle, i.e., with the hyperplane section. 
The Hard Lefschetz theorem, which is valid for intersection
homology by \cite{BBD}, asserts that $\Lambda$ is an isomorphism
in both cases indicated by $\cong$\,.
Thus $p_*$ and $q^*$ vanish, so $p^*$ is
surjective and $q_*$ is injective.
Analyzing the diagram, we conclude that the map
$(\alpha_{\l_X  ,\l_Y})_k$ vanishes.\qed\enddemo

\remark{Remark 3.2} Modifying our method of proof
we can generalize Proposition 3.1 for 
quasihomogeneous singularities.\endremark \vskip5pt

Now we are ready to prove Theorem 1.2.

\demo{Proof of Theorem 1.2} The proof is an induction on $\dim X$.
For $\dim X=0$, the theorem is obvious.

Suppose that $\dim X=n$. We want to construct a morphism
$\alpha_{X,Y}$ which is an extension of $\alpha_\roman{reg}\,$.
We do it stratum by stratum going down with respect to the
partial ordering among the strata. To extend it over
$S$, we choose a slice $N_S$ given by the definition of conical
singularities. 
By the inductive assumption there exists a morphism of intersection
homology $\alpha'_{\goth X,\goth Y}:\i {\goth X}[-2]@>>>\i {\goth Y}$,
where $\goth X=\l_X/S^1$ and $\goth Y=\l_Y/S^1$. By the uniqueness,
the morphism $Rq^*\alpha'_{\goth X,\goth Y}$ coincides with
$\alpha_{\l_X,\l_Y}={\alpha_{X',Y'}}_{|\l_Y}$. (Here $q:\l_Y@>>> \goth
Y$ is again the quotient map.)
By Proposition 3.1, the map $(\alpha_{\l_X \l_Y})_k$ vanishes.
Thus by Proposition 2.2, an extension exists and is unique. Since
$\alpha_\roman{reg}\,$ was geometric, the obtained extension is
geometric as well.
\qed\enddemo

In the end of this section we show how to deduce the existence of a
morphism $\alpha_{X,Y}$ for a pair of analytic singularities. Our
argument is based on the theory of mixed Hodge modules due to
M.~Saito. The proof is modified version of \cite{We}.

\proclaim{Theorem 3.3} Assume that $(Y,X)$ is a pair of analytic
varieties with codimension $\codim X Y=1$ and suppose that no 
component of $X$ is
contained in the singularities of $Y$. Then there exists an unique
extension $\alpha_{X,Y}$ of $\alpha_\roman{reg}$.\endproclaim

\demo{Proof} The question about the existence of $\alpha_{X,Y}$ has
a local nature. 
If it exists locally, then the morphism of links vanishes.
Then, by Proposition 2.2, a global morphism exists, it is unique
and locally it coincides with the original local ones.
Thus we can always patch morphisms. We will construct
$\alpha_{X,Y}$ locally.
Analyzing the techniques of resolution of singularities, we note that
for each $x\in X$ one can find a neighbourhood $U$ in $Y$ and a
resolution $\pi : \tilde U@>>>U$, which is a projective morphism,
i.e.~it factors through an inclusion and the projection
$$ \tilde U \hookrightarrow U\times \Pt^N\epi{} U\,.$$
The decomposition theorem is valid for $\pi$ by 
\cite{S2}, 0.3 or 0.6. Hence $\i U$ is a direct summand in $R\pi_*\D {\tilde
U}$. Analogously $\i {U\cap X}$ is a direct summand in $R\pi_*\i
{\pi^{-1}(U\cap X)}$. Then the composition 
$$\i {U\cap X}[-2]\hookrightarrow 
R\pi_*\i{\pi^{-1}(U\cap X)}[-2]@>R\pi_*\roman{incl}_\#>>
R\pi_*\D {\tilde U}\epi{}\i U$$
is an extension of $\alpha_\roman{reg}$ after suitable rescaling.\qed\enddemo

\remark{Remark 3.4} As before if we allow $X$ to lie in the
singularities of $Y$, then $\alpha_{X,Y}$ exists, but it depends on
the choice made on $X_\roman{reg}$. 
\iii In fact we can proceede as follows.
Suppose $X$ equals the closure of a singular stratum $Y_\alpha\subset Y$.
Let $L_\alpha=\Cal H^0(\i Y)_{|Y_\alpha}$ be 
the cohomology sheaf. It is locally constant, i.e.~it is a
coefficient system on $Y_\alpha$. Then we obtain a morphism
$\i X(L_{\alpha}^*)[-2]@>>>\i Y$, which is an unique extension of the
dual of the restriction morphism on $Y_\roman{reg}\cup Y_\alpha$:
$$\i {Y_\roman{reg}\cup Y_\alpha}@>>>\i {Y_\alpha}(L_\alpha)\cong
L_\alpha\,.$$ \endremark

\head 4. Example -- Projective Cones\endhead

Let ${\goth X}\subset {\goth Y}\subset \Pt^N$ be a pair 
of complex
projective varieties, $\dim {\goth X}=n-1$, $\dim {\goth Y}=n$. Let $K {\goth 
X}\subset
\Pt^{N+1}$ and $K {\goth Y}\subset \Pt^{N+1}$ be the projective cones over
${\goth X}$
and ${\goth Y}$. Topologically they are the Thom spaces of the tautological
bundle restricted to ${\goth X}$ and ${\goth Y}$. 
They are also called algebraic suspensions. 
Suppose we are given a map
$$(\alpha_{{\goth X},{\goth Y}})_*:\h * {\goth X}@>>> \h * {\goth Y}\,.$$
A question arises~: What is the map
$$ (\alpha_{{K\goth X},K{\goth Y}})_*:\h * {K {\goth X}}@>>> \h * {K {\goth 
Y}}\,,$$
which comes from the continuation of $\alpha_{{\goth X},{\goth Y}}$?
We give an answer in the following table. The intersection homology
of a projective cone is computed in \cite{FK}, 3.5. Let $x_i=\h i {\goth X}$ 
and
$y_i=\h i {\goth Y}$.\vskip5pt

\centerline{\vbox{\offinterlineskip
\halign{\strut\vrule\;\hfil $#$\hfil\;&
              \vrule\hfil $#$ \hfil &
              \vrule\;\hfil $#$ \hfil\; &
              \vrule \;\hfill $#$&$#$ \hfill\; & 
              \vrule\;\hfil $#$ \hfil\; & 
              \vrule\hfil $#$ \hfil \vrule \cr
\noalign{\hrule}
\dim & & \h *{K {\goth X}} & & \;\;\roman {map} & \h *{K {\goth Y}} & \cr
\noalign{\hrule\vskip1pt\hrule}
0 & \;* & x_0 & & (\alpha_{{\goth X},{\goth Y}})_0 & y_0 & \; *\cr
\noalign{\hrule}
1 &\; * & x_1 & & (\alpha_{{\goth X},{\goth Y}})_1 & y_1 & \;*\cr
\noalign{\hrule}
\dots & \; * &&&&&\;*\cr\noalign{\hrule}
n-1 &\; * & x_{n-1} & & (\alpha_{{\goth X},{\goth Y}})_{n-1} & y_{n-1} &\; 
*\cr
\noalign{\hrule}
n & \;* & \roman{im}(\Lambda:x_n @>\cong>> x_{n-2})& &( 
\overline{\alpha_{{\goth
X},{\goth Y}}})_n & y_n & \;*\cr
\noalign{\hrule}
n+1 & \;K & x_{n-1} & \Lambda^{-1}&(\alpha_{{\goth X},{\goth Y}})_{n-1} & 
\roman{im}(\Lambda:y_{n+1}@>\cong>> y_{n-1}) &\; 
*\cr
\noalign{\hrule}
n+2 & \;K & x_{n} & & (\alpha_{{\goth X},{\goth Y}})_n & y_{n} &\;K \cr
\noalign{\hrule}
\dots& \;K &&&&&\;K \cr\noalign{\hrule}
2n & \;K & x_{2n-2} & & (\alpha_{{\goth X},{\goth Y}})_{2n-2} & y_{2n-2} &\;K 
\cr
\noalign{\hrule}
2n+1 & & 0 & & & y_{2n-1} &\;K \cr
\noalign{\hrule}
2n+2 & & 0 & & & y_{2n} &\;K \cr
\noalign{\hrule}
}}}\vskip5pt
The star $*$ in the table indicates, that the classes of the
corresponding intersection homology group are represented by the
cycles contained in 
the base of the cone otherwise they are represented by the cones
over the cycles. In the dimension $n$ the map $(\alpha_{{\goth
X},{\goth Y}})_n$ defines a map $(\overline{\alpha_{{\goth X},{\goth 
Y}}})_n:\roman{im}(\Lambda:x_n@>\cong>>x_{n-2})@>>>y_n$. In the
dimension $n+1$ we apply $({\alpha}_{{\goth X},{\goth Y}})_{n-1} $ to
a cycle from 
$x_{n-1}$ and then apply the inverse of $\Lambda$ to obtain an
allowable representative in $y_{n+1}$.

\head 5. Vanishing in Homology of Link Mappings\endhead

Goresky and MacPherson \cite{GM2} list homological properties of
algebraic varieties which follow from the decomposition theorem
of \cite{BBD}. We will add to this list another one.
Proposition 2.2 says that a remarkable property of links is
responsible for the continuation of $\alpha_{X', Y'}$ on $S$.
By Theorem 3.3 and Remark 3.4 such a continuation exists for arbitrary
analytic varieties. 

\proclaim{Corollary 5.1} Let $(Y,X)$ be a pair of analytic
varieties with $\codim X Y=1$. 
Let $\l_X$ and $\l_Y$ be the links of a stratum of codimension $k$ in $X$
and let $\alpha_{X\m \overline S,Y\m \overline S}:\i {X\m \overline
S}[-2]@>>> \i {Y\m \overline S}$ be any morphism of intersection homology. 
Then the induced map 
$(\alpha_{\l_X,\l_Y})_k: \h k {\l_X}@>>>\h k {\l_Y}$ vanishes.
\endproclaim

In the algebraic case the argument of \cite{BBFGK}, p. 167 and Lemme clef
3.3, is to show vanishing of
certain sheaf morphism which directly implies Corollary 5.1.
For a pair of analytic varieties of codimension $d>1$ a morphism
$\alpha_{X,Y}:\i X[-2d]@>>>\i Y$ can be constructed locally as
in the proof of 3.3. (There is no 
guarantee that they come from a global morphism.)
The diagram 2.4 proves vanishing as before. For isolated
singularities, we obtain: 

\proclaim{Corollary 5.2} Let $X\subset Y$ be a pair of analytic
varieties with $\dim X=n$ and $\dim Y=n+d$. Suppose that $x$ is an isolated 
singular point of $X$ and $Y$.
Let $\l_X\subset\l_Y$ be its links.
Then the maps induced by the inclusion 
$\iota_*: H_i(\l_X)@>>>H_i (\l_Y)$ vanish for
$n\leq i < n+d$.
\endproclaim

Corollary 5.2 is related to the results of \cite{HL}, although it is of 
different nature. \iii It can be derived from the fact, that $H^*(\l_X)$
and $H^*(\l_Y)$ are equipped with the mixed Hodge structure, see e.g.
\cite{Di}, C28:
$$H^i(\l_X) \roman{~has~all~the~weights~}\left\{\matrix 
\leq i\roman{~for~} i< n\hfill \\
>i\roman{~for~} i\geq n\,, \endmatrix\right.$$
$$H^i(\l_Y) \roman{~has~all~the~weighs~}\left\{\matrix 
\leq i\roman{~for~} i< n+d \hfill\\
>i\roman{~for~} i\geq n+d\,. \endmatrix\right.$$
The map $\iota^*:H^*(\l_Y)@>>>H^*(\l_X)$ preserves mixed Hodge
structure, thus it has to 
vanish in the range $n\leq i <n+d$, as noticed in the proof of
\cite{Di}, C34. For nonisolated singularities the analogous statement
is also true. By \cite{S1}, 1.18 or \cite{DS} the intersection cohomology of links is
equipped with mixed Hodge structure with the weights as above.
If $(\alpha_{\l_X,\l_Y})^*$ preserves the weights, then it has to
vanish in the middle.

\remark{Remark 5.3} 
If $L$ is a local system defined on an open set $U\subset 
Y_\roman{reg}\,$ and no component of $X$ is contained in $Y\m U$, then
by standard operation on sheaves (as in Section 2)
we obtain a morphism of intersection homology with twisted coefficients:
$$\roman{Incl}_\#:\i X(L_{|X\cap U})[-2]@>>>\p Y(L)\,.$$
Suppose that a local system has geometric origin (\cite{BBD}, p.~162
or \cite{S2}, p.~128)
then modifying the proof of \cite{BBFGK} or \cite{We}
one can show generalized Corollaries 5.1--2 for
twisted coefficients. This way we can prove 
the existence of a morphism 
$$(\alpha_{X,Y})_*:\i X(L_{|X\cap U})[-2]@>>>\i Y(L)\,.$$
\endremark\vskip5pt

The Corollaries 5.1--2 we consider as a ``local Hard
Lefschetz'' theorem. The reason for this is following:

\proclaim{Proposition 5.4} The Corollary 5.1 implies the Hard Lefschetz
theorem for intersection homology.\endproclaim

\demo{Proof} Let $\goth Y$ be a projective variety of dimension $n$
and let $\goth X=\goth Y\cap H$ be a generic hyperplane section. The inclusion
$i:\goth X\hookrightarrow \goth Y$ is normally nonsingular, thus it
induces a map of intersection homology. Consider again the diagram with Gysin
sequences of 3.1:
$$\CD @>>>\h{n+1}{\goth X}@>\Lambda>>\h{n-1}{\goth X}@>p^*>>\h n{\l_X}@>>>\\
@. @VVi_{n+1}V @VVi_{n-1}V @VV0V\\
@>>>\h{n+1}{\goth Y}@>\Lambda>>\h{n-1}{\goth Y}@>q^*>>\h 
n{\l_Y}@>>>\,.\endCD$$
By the weak Lefschetz theorem
\cite{GM3}, \S6.10 the map $i_{n-1}$ is surjective. Thus if the map of
links vanishes in the middle homology, then $q^*$ vanishes as well.
We conclude that $\Lambda$ for $\goth Y$ is surjective. By Poincar\'e duality 
it
is an isomorphism. The argument that $\Lambda^k:\h{n+k}{\goth
Y}@>>>\h{n-k}{\goth X}$ is an isomorphism for $k>1$ is standard; it
follows from the weak Lefschetz theorem.\qed\enddemo

\head 6. Application -- Chern Classes\endhead

Let $X$ be a $n$--dimensional analytic variety. Suppose we are given a
sequence of subvarieties:
$$ \cF =\{X=X^0\supset X^1\supset\dots\supset X^k\}\,. $$
Assume that no component of $X^i$ is contained in the singularities
of $X^{i-1}$ and $\codim  {X^i} X =i$ for
$i=1,2,\dots , k$. 
We will call such sequence {\it a flag} in $X$. With a flag in $X$, we
associate a sequence of  
elements in intersection homology:
$[\cF ]^i\in \hc {2(n-i)}(X)$. 
These classes are lifts of the corresponding fundamental classes
$[X^i] \in H^\roman{cld}_{2(n-i)}(X)$. The construction of $[\cF]^i$
is the following~: We take the fundamental class of $X^i$ in $\hc
{2(n-i)}(X^i)$ and we lift it step by step 
going through $\hc {2(n-i)}(X^j)$ for $0\leq j<i$.
On each step we choose the lift induced by
$(\alpha_{X^{j},X^{j-1}})_*$.

\remark{Example 6.1} 
Let $\goth X\subset \Pt^N$ be a projective variety of pure dimension
$n$ and let $c\l_X\subset\Ct^{N+1}$ be its affine cone.
Consider a flag of linear subspaces
$$\Cal V=\{V^{n+1}\subset V^n\subset \dots\subset V^1\subset
V^0=\Ct^{N+1}\}$$ 
with $\codim {V^i}{\Ct^{N+1}}=i$ for $i=0,1,\dots,n+1$.
Let $$p_i:c\l_X@>>>\Ct^{N+1}/V^{n+2-i}\cong\Ct^{n+2-i}$$ 
be the linear projection restricted to $c\l_X$. Consider the
set of critical points of $p_i$
$$(c\l_X)^i=\roman{closure}{\left\{x\in(c\l_X)_\roman{reg}: 
 p_i\;\roman{is\; not\; submersion\; at}\;x\right\}}\,,$$
where $(c\l_X)_\roman{reg}\,$ is the nonsingular part of $c\l_X$.
Each set $(c\l_X)^i$ is homogeneous.  We denote its projectivization
by $\goth X^i$. We call it the $i$-th polar variety of $\goth X$. We have
$\goth X^i\subset \goth X^{i-1}$ for $0<i\leq n+1$ and we put $\goth
X^0=\goth X$. If the flag $\Cal V$ is general enough, then the
sequence $$\cP=\{\goth X^0=\goth X\supset\goth
X^1\supset\dots\supset\goth X^n\}$$ is a flag in $\goth X$ and $\goth
X^{n+1}=\emptyset$.  We obtain classes $[\cP]^i\in \h{2(n-i)}{\goth
X}$. If $\goth X$ is smooth, then $(-1)^i[\cP]^i=(-1)^i[\goth X^i]\in
\h{2(n-i)}{\goth X}
\cong H_{2(n-i)}(\goth X)$ is the Poincar\'e dual of the Chern class 
of a bundle $\Theta$,
where $\Theta$ fits to an exact sequence of bundles:
$$0@>>>\Cal O(-1)@>>>\Theta@>>>T\goth X\otimes\Cal O(-1)@>>>0\,,$$
see \cite{Po}. Then the duals of the Chern classes of $\goth X$
satisfy the formula: 
$$c_i(\goth X)=\sum_{j=0}^i \binom{n+1-j}{i-j}(-1)^jh^{i-j}\cap
[\goth X^j]\,,\tag 6.2$$ where $h\in H^2(\goth X)$ is the class of
the hyperplane section,
compare \cite{Fu}, \S14.4.15. When $\goth X$ is singular, then by \cite{Pi}
the formula 6.1 describes the Chern-Mather class~:
$c^M_i(\goth X)\in H_{2(n-i)}(\goth X)$ 
which in general does not lie in the image of the Poincar\'e duality map.
By the very same formula with $[\goth X^j]$ replaced by $[\cP]^j$ we
define an element in intersection homology. This is
a particular lift of $c_i^M(\goth X)$ to intersection homology. We
conjecture that it does not  
depend on the embedding.\endremark

\remark{Remark 6.3} This method of defining characteristic classes was already
known in the thirties. This  is  the  way  how  Todd in  \cite{Td} defined  a 
"canonical system" with the hyperplane bundle replaced by an
arbitrary bundle. It was later called the Chern class of $\goth X$.
There are enough evidences that polar varieties carry an important
information about global and local invariants of $\goth X$, see 
e.g.~\cite{LT}.
The Chern classes in intersection homology were defined by J-P.~Brasselet and 
by
the author. The paper \cite{BW} contains a proof that they do not
depend on the choice of the generic flag $\Cal V$. It also
contains an explicit computation for $\goth X$ which is the
projective cone over a 
quadric in $\Pt^3$.\endremark

\enddocument
\end